\documentclass[10pt]{article}

\usepackage{a4wide}
\usepackage{amssymb}
\usepackage{amsfonts}
\usepackage{amsmath}
\input xy
\xyoption{arrow} \xyoption{matrix}

\date{}

\newtheorem{proposition}{Proposition}[section]
\newtheorem{theorem}[proposition]{Theorem}
\newtheorem{lemma}[proposition]{Lemma}

\newtheorem{corollary}[proposition]{Corollary}

\def\Hom{{\rm Hom}}
\def\der{\partial }

\def\nFM0{{\nu }_{F,M_0}}
\def\nFN0{{\nu }_{F,N_0}}
\def\nGN0{{\nu }_{G,N_0}}

\def\N0{ {\bf N}_0 }

\def\ra{\rightarrow}

\def\Xpm{X^{\pm }}

\def\s{\sigma}
\def\Z{\mathbb{Z}}

\def\l1{{\lambda}_1}

\def\a{\alpha}
\def\a0{ {\alpha }_0}
\def\a1{ {\alpha }_1}

\def\l{\lambda}

\def\nFGM0{{\nu }_{F,G,M_0}}

\def\nFN0{{\nu}_{F,N_0}}

\def\sm{{\sigma}^m}

\def\sm1{{\sigma}^{-1}}

\def\smtp1{{\sigma}^{-t+1}}

\def\S1{S^{-1}}

\def\Xpm1{X^{\pm 1}_1}

\def\sPM1{{\sigma }^{\pm 1}}
\def\sMP1{{\sigma }^{\mp 1 }}

\def\d{\delta}

\def\di{{\rm d.ind}}

\def\L{\Lambda}

\def\CD{{\cal D}}

\def\Ytm1{Y^{t-1}}
\def\Yim1{Y^{i-1}}

\def\CK{{\cal K}}

\def\CF{{\cal F}}
\def\CG{{\cal G}}
\def\CH{{\cal H}}

\def\Aut{{\rm Aut}}

\def\Der{{\rm Der }}
\def\ad{{\rm ad }}
\def\dim{{\rm dim }}

\def\ker{ {\rm ker } }

\def\SL2Z{ {\rm SL}_2({\bf Z}) }

\def\th{ \theta }

\def\Gp1{ G^{1 , 1 } }
\def\P11{ P^{-1 , 1 } }
\def\Pp1{ P^{1 , 1 } }

\def\Autgr{{\rm Aut }_{gr}}

\def\NP{{\rm NP}}
\def\Supp{{\rm Supp}}

\def\th{\theta}

\def\nCLsr{{}^\nu\kern-2pt {\cal L}^{\sigma , \rho  }}
\def\nP{{}^\nu \kern-2pt P}
\def\nL{{}^\nu\kern-2pt L}
\def\nLL{{}^\nu\kern-2pt \Lambda}
\def\nPsr{{}^\nu\kern-2pt P^{\sigma , \rho  }}
\def\nLsr{{}^\nu\kern-2pt L^{\sigma , \rho  }}
\def\nuCL{{}^\nu\kern-2pt  {\cal L}}
\def\nCLsr{{}^\nu\kern-2pt {\cal L}^{\sigma , \rho  }}
\def\nCL1m{{}^\nu\kern-2pt {\cal L}^{-1 , 1  }}
\def\x1nu{x^\frac{1}{\nu}}
\def\xm1nu{x^{-\frac{1}{\nu}}}

\def\ra{\rightarrow }

\def\CB{{\cal B}}

\def\CH{ {\cal H}}

\def\nAM0{{\nu }_{{\cal A},M_0}}
\def\nAN0{{\nu }_{{\cal A},N_0}}

\def\End{ {\rm End }}
\def\Der{ {\rm Der }}

\def\ad{ {\rm ad }}

\def\GL{{\rm GL}}
\def\SL{{\rm SL}}

\def\Hom{{\rm Hom}}

\def\di!{\frac{\der^i}{i!}}
\def\dik!{\frac{\der^k_i}{k!}}

\def\id{{\rm id}}

\def\N{\mathbb{N}}

\def\0{\overline{0}}
\def\1{\overline{1}}

\def\Ln1{\L_{n,\overline{1}}}

\def\a1{a_{\overline{1}}}

\def\S{\Sigma}

\def\vn1{\overrightarrow{n-1}}

\def\Sh{{\rm Sh}}

\def\im{{\rm im}}

\def\mL{\mathbb{L}}

\def\mJ{\mathbb{J}}
\def\mI{\mathbb{I}}

\def\mT{\mathbb{T}}

\def\mG{\mathbb{G}}

\def\mX{\mathbb{X}}

\def\K1{{\rm K}_1}

\def\hmI1{\widehat{\mI_1}}
\def\tmI1{\widetilde{\mI_1}}
\def\tmJ1{\widetilde{\mJ_1}}
\def\hB1{\widehat{B_1}}
\def\hCB1{\widehat{\CB_1}}

\def\ggu{\mathfrak{u}}

\def\Fix{{\rm Fix}}

\def\mJ{\mathbb{J}}

\def\Nil{{\rm Nil}}

\def\divn0{\mathfrak{div}_n^0}
\def\div0mu{\mathfrak{div}_{n, [\mu ]}^0}
\def\din0{\mathfrak{di}_n^0}
\def\ivn0{\mathfrak{iv}_n^0}

\def\AutLie{ {\rm Aut}_{ {\rm Lie} } }
\def\AutKalg{ {\rm Aut_{K-{\rm alg}}}}
\def\Vir{{\rm Vir}}
\def\mmWn{\mathbb{W}_n}
\def\mmW1{\mathbb{W}_1}
\def\mLn{\mL_n}
\def\mX{\mathbb{X}}
\def\LF{{\rm LF}}
\def\diag{{\rm diag}}
\def\Nil{{\rm Nil}}

\begin{document}

\author{V. V. \  Bavula   
}

\title{The groups of automorphisms of the Witt $W_n$ and Virasoro Lie algebras }

\maketitle

\begin{abstract}
Let $L_n=K[x_1^{\pm 1} , \ldots , x_n^{\pm 1}]$ be a Laurent polynomial algebra over a field $K$ of characteristic zero,  $W_n:= \Der_K(L_n)$, the {\em Witt} Lie algebra, and $\Vir$ be the {\em Virasoro} Lie algebra. We prove that $\AutLie (W_n) \simeq \AutKalg (L_n)\simeq \GL_n(\Z ) \ltimes K^{*n}$ and $\AutLie (\Vir ) \simeq \AutLie (W_1) \simeq \{\pm 1\} \ltimes K^*$.

$\noindent $



\end{abstract}


\section{Introduction}

In this paper, module means a left module, $K$ is a
field of characteristic zero and  $K^*$ is its group of units, and the following notation is fixed:
\begin{itemize}
\item $P_n:= K[x_1, \ldots , x_n]=\bigoplus_{\alpha \in \N^n}
Kx^{\alpha}$ is a polynomial algebra over $K$ where
$x^{\alpha}:=x_1^{\alpha_1}\cdots x_n^{\alpha_n}$,
 \item $G_n:=\AutKalg (P_n)$ is the group of automorphisms of the polynomial algebra $P_n$,
 \item $L_n:=K[x_1^{\pm 1} , \ldots , x_n^{\pm 1}]=\bigoplus_{\alpha \in \Z^n} Kx^\alpha $ is a Laurent polynomial algebra,
 \item $\mL_n:=\AutKalg (L_n)$ is the group of $K$-algebra automorphisms of $L_n$,
 \item $\der_1:=\frac{\der}{\der x_1}, \ldots , \der_n:=\frac{\der}{\der
x_n}$ are the partial derivatives ($K$-linear derivations) of
$P_n$,
\item    $D_n:=\Der_K(P_n) =\bigoplus_{i=1}^nP_n\der_i$ is the Lie
algebra of $K$-derivations of $P_n$ where $[\der , \d ]:= \der \d -\d \der $,
\item $\mG_n:=\Aut_{{\rm Lie}}(D_n)$ is the group of automorphisms of the Lie algebra $D_n$,

    \item    $W_n:=\Der_K(L_n) =\bigoplus_{i=1}^nL_n\der_i$ is the {\em Witt}  Lie
algebra where $[\der , \d ]:= \der \d -\d \der $,
\item $\mmWn:=\AutLie (W_n)$ is the group of automorphisms  of the Witt Lie algebra $W_n$,
\item  $\d_1:=\ad (\der_1), \ldots , \d_n:=\ad (\der_n)$ are the inner derivations of the Lie algebras $D_n$ and $W_n$  determined by  $\der_1, \ldots , \der_n$ where $\ad (a)(b):=[a,b]$,

 \item $\CD_n:=\bigoplus_{i=1}^n K\der_i$,
 \item $\CH_n :=\bigoplus_{i=1}^n KH_i$ where $H_1:=x_1\der_1, \ldots , H_n:=x_n\der_n$,


    \end{itemize}

{\bf The group of automorphisms of the Witt Lie algebra $\mmWn$}.
The aim of the paper is to find the groups of automorphisms of the Witt algebra $W_n$ (Theorem \ref{20Mar13}) and the Virasoro algebra $\Vir$ (Theorem \ref{A20Mar13}).  The following lemma is an easy exercise.
\begin{itemize}
\item {\rm (Lemma \ref{g21Mar13})} $\mLn \simeq \GL_n(\Z ) \ltimes \mT^n$ {\em where $\GL_n(\Z )$ is identified with a subgroup of $\mLn$ via the group monomorphism} $ \GL_n(\Z ) \ra \mLn$, $ A= (a_{ij})\mapsto \s_a: x_i\mapsto \prod_{j=1}^nx_j^{a_{ji}}$ {\em and} $$\mT^n:=\{ t_\l \in \mLn\, | \, t_\l (x_1)=\l_1x_1, \ldots , t_\l (x_n) =\l_n  x_n; \l\in K^{*n}\}\simeq K^{*n}$$ {\em is the algebraic $n$-dimensional torus}.
\end{itemize}

\begin{theorem}\label{20Mar13}
$\mmWn = \mLn$.
\end{theorem}
{\em Structure of the proof}.   (i) $\mLn$ is a subgroup $\mmWn$ (Lemma \ref{a21Mar13})
 via the group monomorphism
$$\mLn\ra \mmWn, \;\;  \s \mapsto \s : \der \mapsto \s (\der ):=  \s \der \s^{-1}.$$
Let $\s \in \mmWn$. We have to show that $\s \in \mLn$.

(ii)(crux)  $\s (\CH_n) = \CH_n$ (Lemma \ref{d21Mar13}), i.e.
$$ \s (H) = A_\s H \;\; {\rm for\; some }\;\; A_\s \in \GL_n(K)$$ where $H:=(H_1, \ldots , H_n)^T$.

$\noindent $

(iii) $A_\s \in \GL_n(\Z )$  (Corollary \ref{f21Mar13}).

$\noindent $

(iv) There exists an automorphism $\tau \in \mLn$ such that $\tau \s \in \Fix_{\mmWn}(H_1, \ldots , H_n)$ (Lemma  \ref{k21Mar13}).

$\noindent $

(v) $\Fix_{\mmWn}(H_1, \ldots , H_n)=\mT^n \subseteq \mLn$ (Lemma  \ref{l21Mar13}) and so $\s \in \mLn$. $\Box $

$\noindent $

{\bf The group of automorphisms of the Virasoro Lie algebra}. The {\em Virasoro} Lie algebra $\Vir = W_1 \oplus Kc$  is a 1-dimensional central extension of the Witt Lie algebra $W_1$ where $Z(\Vir ) = Kc$ is the centre of $\Vir$ and for all $i,j\in \Z$,
\begin{equation}\label{defVir}
[x^iH, x^jH]=(j-i)x^{i+j}H+\d_{i,-j} \frac{i^3-i}{12}c
\end{equation}
where $x=x_1$ and  $H=H_1$.

\begin{theorem}\label{A20Mar13}
$\AutLie (\Vir ) \simeq \mmW1\simeq \mL_1\simeq \GL_1(\Z ) \ltimes \mT^1$.
\end{theorem}

The key point in the proof of Theorem \ref{A20Mar13}  is to use Theorem \ref{22Mar13} of which Theorem \ref{A20Mar13}  is a special case (where $\CG = \Vir$, $W= W_1$ and $Z= Kc$,  see Section \ref{AUT-Vir}).

\begin{theorem}\label{22Mar13}
Let $\CG$ be a Lie algebra, $Z$ be a subspace of the  centre  of $\CG$ and $W=\CG / Z$. Suppose that
\begin{enumerate}
\item every automorphism $\s$ of the Lie algebra $W$ can be extended to an automorphism $\widehat{\s}$ of the Lie algebra $\CG$,
\item $Z\subseteq [G, G]$, and
\item $W= [W, W]$.
\end{enumerate}
Then, for each $\s$,  the extension $\widehat{\s}$ is unique  and the map $ \AutLie (W)\ra \AutLie (\CG )$, $\s \mapsto \widehat{\s}$, is a group isomorphism.
\end{theorem}


The groups $\AutLie (\ggu_n)$ and $\AutLie (D_n)$ where found in \cite{Bav-Lie-Un-AUT} and \cite{Bav-Aut-Der-Pol} respectively. The Lie algebras  $\ggu_n$ have been studied in great detail in \cite{Bav-Lie-Un-MON} and \cite{Bav-Lie-Un-GEN}. In particular, in \cite{Bav-Lie-Un-MON} it was proved that every monomorphism of the Lie algebra $\ggu_n$ is an automorphism  but this is not true for epimorphisms.



\section{Proof of Theorem \ref{20Mar13} }\label{P20AAA}

This section can be seen as a proof of Theorem \ref{20Mar13}. The proof is split into several statements that reflect `Structure of the proof of Theorem \ref{20Mar13}' given in the Introduction.

By the very definition, $\CH_n = \bigoplus_{i=1}^n KH_i$ is an abelian Lie subalgebra of $W_n$ of dimension $n$. Each element $H$ of $\CH_n$ is a unique sum $H = \sum_{i=1}^n \l_iH_i$ where $\l_i\in K$. Let us define the bilinear map
$$\CH_n\times \Z^n \ra K, \;\; (H, \alpha ) \mapsto (H, \alpha ) := \sum_{i=1}^n \l_i \alpha_i.$$

{\bf The Witt algebra $W_n$ is a $\Z^n$-graded Lie algebra}. The Witt  algebra
\begin{equation}\label{HbHa}
W_n =\bigoplus_{\alpha\in \Z^n} \bigoplus_{i=1}^n Kx^\alpha \der_i=\bigoplus_{\alpha\in \Z^n}x^\alpha\CH_n
\end{equation}
is a $\Z^n$-graded Lie algebra, that is $[x^\alpha \CH_n , x^\beta \CH_n]\subseteq x^{\alpha + \beta}\CH_n$ for all $\alpha ,\beta \in \Z^n$. This follows from the identity
\begin{equation}\label{HbHa1}
[x^\alpha H, x^\beta H']= x^{\alpha + \beta } ((H,\beta )H'-(H',\alpha ) H).
\end{equation}
In particular,
\begin{equation}\label{HxaH}
[ H, x^\alpha H']= (H, \alpha ) x^\alpha H'.
\end{equation}
So, $x^\alpha \CH_n$ is the {\em weight subspace} $W_{n,\alpha } := \{ w\in W_n \, | \, [ H, w] = (H, \alpha ) w\}$  of $W_n$ with respect to the adjoint action of the abelian Lie algebra $\CH_n$ on $W_n$. The direct sum (\ref{HbHa}) is the weight decomposition of $W_n$ and $\Z^n$ is the {\em set of weights} of $\CH_n$.

Let $\CG$ be a Lie algebra and $\CH$ be its Lie subalgebra. The {\em centralizer} $C_\CG (\CH ) := \{ x\in \CG \, | \, [ x, \CH ] =0\}$ of $\CH$ in $\CG$ is a Lie subalgebra of $\CG$. In particular, $Z(\CG ) := C_{\CG }(\CG ) $ is the {\em centre} of the Lie algebra $\CG$. The {\em normalizer} $N_\CG (\CH ) :=\{ x\in \CG \, | \, [ x, \CH ] \subseteq \CH\}$ of $\CH$ in $\CG$ is a Lie subalgebra of $\CG$, it is the largest Lie subalgebra of $\CG$ that contains $\CH $ as an ideal.  Each element $a\in \CG$ determines the
derivation  of the Lie algebra $\CG$ by the rule $\ad (a) : \CG
\ra \CG$, $b\mapsto [a,b]$, which is called the {\em inner
derivation} associated with $a$. An element $a\in \CG$ is called a {\em locally finite element}  if so is the inner derivation $\ad (a)$ of the Lie algebra $\CG$, that is $\dim_K(\sum_{i\in \N} K\ad (a)^i (b))<\infty$ for all $b\in \CG$. Let $\LF (\CG )$ be the set of locally finite elements of $\CG$.

$\noindent $

{\bf The Cartan subalgebra $\CH_n$ of $W_n$}. A nilpotent Lie subalgebra $C$ of a Lie algebra $\CG$ such that $ C=N_\CG (C)$  is called a {\em Cartan subalgebra} of $\CG$. We use often the following obvious observation: {\em An abelian Lie subalgebra that coincides with its centralizer is a maximal abelian Lie subalgebra}.

\begin{lemma}\label{a31Mar13}
\begin{enumerate}
\item $\CH_n=C_{W_n}(\CH_n)$  is a maximal abelian Lie subalgebra of $W_n$.
\item $\CH_n$  is a Cartan subalgebra of $W_n$.
\end{enumerate}
\end{lemma}

{\it Proof}. Both statements follow from (\ref{HbHa}) and (\ref{HxaH}) . $\Box$

The next lemma is very useful and can be applied in many different situations. It allows one  to see the group of automorphisms  of a ring as a subgroup of the group of automorphisms of its Lie algebra of derivations.

\begin{lemma}\label{a21Mar13}
Let $R$ be a commutative ring such that there exists a derivation $\der \in \Der (R)$ such that $r\der \neq 0$ for all nonzero elements $r\in R$ (eg, $R= P_n, L_n$ and $\d = \der_1$). Then the group homomorphism
$$ \Aut (R) \ra \AutLie (\Der (R)), \;\; \s \mapsto \s : \d\mapsto \s (\d ) := \s \d \s^{-1},$$ is a monomorphism.
\end{lemma}

{\it Proof}. If an automorphism $\s \in \Aut (R)$ belongs to the kernel of the group homomorphism $\s \mapsto \s$ then,  for all $r\in R$, $ r\der = \s(r\der )\s^{-1} = \s (r)\s\der \s^{-1} = \s (r) \der$, i.e. $\s (r)= r$ for all $r\in R$. This means that $\s $ is the identity automorphism. Therefore, the homomorphism $\s\mapsto \s$ is a monomorphism. $\Box $

$\noindent $

{\bf The $(\Z , \l )$-grading and the filtration  $\CF_\l$ on $W_n$}. Each vector $\l = (\l_1, \ldots , \l_n) \in \Z^n$ determines the $\Z$-grading on the Lie algebra $W_n$ by the rule
$$W_n=\bigoplus_{i\in \Z} W_{n,i}(\l ), \;\; W_{n,i}(\l ) :=\bigoplus_{(\l , \alpha )=i} x^\alpha \CH_n , \;\; (\l , \alpha ) := \sum_{i=1}^n \l_i \alpha_i, $$
 $[W_{n,i}(\l ) , W_{n,j}(\l ) ] \subseteq W_{n, i+j}(\l )$ for all $i,j\in \Z$ as follows from  (\ref{HbHa1}) and (\ref{HxaH}). The $\Z$-grading above is called the $(\Z , \l)$-{\em grading} on $W_n$. Every element $a\in W_n$ is the unique sum of homogeneous elements with respect to the $(\Z , \l )$-grading on $W_n$,
$$ a= a_{i_1}+a_{i_2}+\cdots  + a_{i_s}, \;\; a_{i_\nu } \in W_{n, i_\nu}(\l),$$
and $i_1<i_2<\cdots <i_s$. The elements $l_\l^+(a):= a_{i_s}$  and  $l_\l^-(a):= a_{i_1}$ are called the {\em leading term} and   the {\em least term}  of $a$ respectively. So,
\begin{eqnarray*}
 a&=& l_\l^+(a)+\cdots ,\\
 a&=& l_\l^-(a)+\cdots ,
\end{eqnarray*}
where the three dots denote smaller and larger terms respectively. For all $a,b\in W_n$,
\begin{equation}\label{ablp}
[a,b]=[l_\l^+(a), l_\l^+(b)]+\cdots ,
\end{equation}
\begin{equation}\label{ablm}
[a,b]=[l_\l^-(a), l_\l^-(b)]+\cdots ,
\end{equation}
where the three dots denote smaller and larger terms respectively (the brackets on the RHS can be zero).

$\noindent $

{\bf The Newton polygon of an element of $W_n$}. Each element $a\in W_n$ is the unique finite sum $a=\sum_{\alpha \in \Z^n} \l_\alpha x^\alpha H_\alpha$ were $\l_\alpha \in K$ and $H_\alpha \in \CH_n$. The set $\Supp (a) := \{  \alpha \in \Z^n\, | \, \l_\alpha \neq 0\}$ is called the {\em support} of $a$ and its convex hull in $\mathbb{R}^n$ is called the {\em Newton polygon} of $a$, denoted by $\NP (a)$.

\begin{lemma}\label{c21Mar13}
Let $a$ be a locally finite element of $W_n$. Then the elements $l_\l^+(a)$ and $l_\l^-(a)$ are locally finite for all $\l \in \Z^n$.
\end{lemma}

{\it Proof}. The statement follows from (\ref{ablp}) and (\ref{ablm}).  $\Box $

$\noindent $

Let $\LF (W_n)_h$ be the set of {\em homogeneous} (with respect to the $\Z^n$-grading on $W_n$) locally finite elements of the Lie algebra $W_n$.

\begin{lemma}\label{b21Mar13}
$\LF (W_n)_h=\CH_n$.
\end{lemma}

{\it Proof}. $\CH_n\subseteq \LF (W_n)_h$ since every element of $\CH_n$ is a semi-simple element of $W_n$: for all $H=\sum_{i=1}^n \l_iH_i$ where $\l_i\in K$,
\begin{equation}\label{HxaHp}
[H, x^\alpha H']= (\l , \alpha ) x^\alpha H' \;\; {\rm for\; all}\;\; \alpha \in \Z^n, \; H'\in \CH_n.
\end{equation}
It suffices to show that every homogeneous element $x^\alpha H'$ that does not belong to $\CH_n$, i.e. $\alpha\neq 0$, is not locally finite. Fix $i$ such that $\alpha_i\neq 0$. Let $\d = \ad (x^\alpha H')$.

Suppose that $(H',\alpha )\neq 0$. This is the case for $n=1$. Then
$$\d^m (x^{2\alpha }H') = (m-1)!2^{m-1} (H', \alpha )^m x^{(1+2m)\alpha } H'\;\; {\rm for}\;\; m\geq 1. $$
Therefore, the element $x^\alpha H'$ is not locally finite.

Suppose that $(H', \alpha ) =0$. Then necessarily $n\geq 2$. Fix $\beta \in \Z^n$ such that $(H' , \beta ) =1$. Then
$$ \d^m (x^\beta H') = x^{\beta +m\alpha } H'\;\; {\rm for}\;\; m\geq 1.$$
Therefore, the element $x^\alpha H'$ is not locally finite. $\Box $


\begin{lemma}\label{d21Mar13}
$\s (\CH_n) = \CH_n$ for all $\s \in \mmWn$.
\end{lemma}

{\it Proof}. Let $\s \in \mmWn$ and $H\in \CH_n$. We have to show that $H':= \s (H) \in \CH_n$. The element $H$ is a locally finite element, hence so is $H'$. By Lemma \ref{c21Mar13} and Lemma \ref{b21Mar13}, the Newton polygon $\NP (H')$ has the single vertex $0$, i.e. $H'\in \CH_n$.  $\Box $

$\noindent $

Let $H= (H_1, \ldots , H_n)^T$ where $T$ stands for the transposition. By Lemma \ref{d21Mar13},

\begin{equation}\label{sHAs}
\s (H) = A_\s H\;\; {\rm for \; all}\;\; \s \in \mmWn
\end{equation}
where $A_\s = (a_{ij})\in \GL_n(K)$ and  $\s (H_i) = \sum_{j=1}^n a_{ij}H_j$. Let ${}^\s W_n$ be the $W_n$-module $W_n$ twisted by the automorphism $\s \in \mmWn$. As a vector space, ${}^\s W_n=W_n$,  but the adjoint action is twisted by $\s $:
$$w\cdot x^\alpha H''= [ \s (w), x^\alpha H'']$$ for all $w\in W_n$ and  $\alpha \in \Z^n$.  The map $\s : W_n\ra {}^\s W_n$, $w\mapsto \s (w)$, is a $W_n$-module isomorphism. By Lemma \ref{d21Mar13}, every  weight subspace $x^\alpha \CH_n$ of the $\CH_n$-module $W_n = \bigoplus_{\alpha \in \Z^n} x^\alpha \CH_n$ is also a weight subspace for the $\CH_n$-module ${}^\s W_n$, and vice versa. Moreover,
\begin{equation}\label{sWnA}
W_{n,\alpha } = x^\alpha \CH_n = ({}^\s W_n)_{A_\s \alpha}\;\; {\rm for \; all}\;\; \alpha \in \Z^n
\end{equation}
where $\alpha = (\alpha_1, \ldots , \alpha_n)^T\in \Z^n$ is a  column: for all $H'= \sum_{i=1}^n \l_i H_i\in \CH_n$,
\begin{equation}\label{sWaA1}
 [ \s (H') , x^\alpha H'']= \sum_{i,j=1}^n \l_i a_{ij} \alpha_j  x^\alpha H''= (H', A_\s \alpha ) x^\alpha H''.
\end{equation}
 Since $\s ( \CH_n ) = \CH_n$ and $\s : W_n\ra {}^\s W_n$ is a $W_n$-module isomorphism, the automorphism $\s$  permutes the weight components $\{ W_{n,\alpha}=x^\alpha \CH_n \}_{\alpha \in \Z^n}$. There is a bijection $\s' : \Z^n\ra \Z^n$, $\alpha \mapsto \s' (\alpha )$, such that $\s (W_{n,\alpha }) = W_{n, \s ' (\alpha )}$ for all $\alpha \in \Z^n$.

\begin{lemma}\label{e21Mar13}
For all $\s \in \mmWn$ and $\alpha \in \Z^n$, $\s'(\alpha ) = A_{\s^{-1}} \alpha $.
\end{lemma}

{\it Proof}.  By (\ref{sWaA1}),
\begin{eqnarray*}
 (H', \s'(\alpha )) \s (x^\alpha H'') &=& [H', \s (x^\alpha H'')]= \s ([\s^{-1}(H'), x^\alpha H''])= \s ((H', A_{\s^{-1}}\alpha) x^\alpha H'')\\
 &=& (H'', A_{\s^{-1}}\alpha ) \s (x^\alpha H'').
\end{eqnarray*}


Therefore, $\s' (\alpha ) = A_{\s^{-1}} \alpha$. $\Box $


\begin{corollary}\label{f21Mar13}
For all $\s \in \mmWn$, $A_\s \in \GL_n(\Z )$.
\end{corollary}

{\it Proof}. This follows from Lemma \ref{e21Mar13}.  $\Box $

$\noindent $

{\bf The group of automorphisms $\mLn = \AutLie (L_n)$}. The group $\mLn$ contains two obvious subgroups: the {\em algebraic $n$-dimensional torus} $\mT^n = \{ t_\l \, | \, \l \in K^{*n}\} \simeq K^{*n}$ where $t_\l (x_i) = \l_ix_i$ for $i=1, \ldots , n$ and $\GL_n(\Z )$ which can be seen as a subgroup of $\mLn$ via the group monomorphism
\begin{equation}\label{msAG}
\GL_n(\Z )\ra \mLn ,\;\;  A\mapsto \s _A: x_i\mapsto \prod_{j=1}^n x_j^{a_{ji}}.
\end{equation}
For all $\alpha \in \Z^n$, $\s_A(x^\alpha ) = x^{A\alpha }$. Hence $\s_{AB}= \s _A\s_B$ and $\s_A^{-1} = \s_{A^{-1}}$.
\begin{lemma}\label{g21Mar13}
$\mLn = \GL_n (\Z ) \ltimes \mT^n$.
\end{lemma}

{\it Proof}. The group of units $L_n^*$ of the algebra $L_n$ is equal to the direct product of its two subgroups $K^*\times \mX$
 where $\mX = \{ x^\alpha \, | \, \alpha \in \Z^n \} \simeq \Z^n$ via $x^\alpha \mapsto  \alpha$. Since $\s (K^*)= K^*$ for all $\s\in \mLn$, there is a group homomorphism (where $\Autgr (G)$ is the group of automorphisms of a {\em group} $G$)
 $$\th :  \mLn \ra \Autgr (L_n/ K^*) , \;\; \s \mapsto \overline{\s} :K^*x^\alpha \mapsto K^*\s (x^\alpha ).$$
 Notice that $\Autgr (L_n/ K^*)\simeq \Autgr (\Z^n)\simeq \GL_n(\Z )$ and $ \th|_{\GL_n(\Z )}: \GL_n(\Z ) \ra \Autgr (L_n/ K^*)$, $A\mapsto A$. Then $\mLn\simeq \GL_n(\Z )\ltimes \ker (\th )$ but $\ker (\th ) = \mT^n$. Clearly, $\mLn = \GL_n (\Z ) \ltimes \mT^n$.  $\Box $


\begin{lemma}\label{h21Mar13}
Let $\s_A\in \mLn$ be as in (\ref{msAG}) where $A\in \GL_n(\Z )$, $\der = (\der_1, \ldots , \der_n)^T$, $H= (H_1, \ldots , H_n)^T$ and   $\diag (\l_{1 1}, \ldots , \l_{nn})$ be  the diagonal matrix with the diagonal elements $\l_{11}, \ldots , \l_{nn}$. Then
\begin{enumerate}
\item $\s (\der ) = C_\s \der$ where $C_\s = \diag (\s (x_1)^{-1},  \ldots , \s(x_n)^{-1})A^{-1} \diag (x_1,  \ldots , x_n)$.
\item $\s (H) = A^{-1} H$.
\end{enumerate}
\end{lemma}

{\it Proof}. 1. Let $\der_i'=\s (\der_i)$ and  $x_j'=\s (x_j)$.  Clearly, $ \s (\der ) = C_\s \der $ for some matrix   $C_\s = (c_{ij})\in M_n(L_n)$. Applying the automorphism $\s$ to the equalities $\d_{ij} = \der_i*x_j$ where  $i,j=1, \ldots , n$, we obtain the equalities $$\d_{ij} = \s\der_i\s^{-1} \s (x_j) = \der_i'*x_j'=(\sum_{k=1}^n c_{ik}\der_k) *\prod_{l=1}^n x_l^{a_{lj}}= (\sum_{k,l=1}^n c_{ik} x_k^{-1}a_{kj})x_j'$$
 where $i,j=1, \ldots , n$. Equivalently, $C_\s \diag (x_1^{-1}, \ldots , x_n^{-1}) A= \diag (x_1'^{-1},\ldots , x_n'^{-1})$, and statement 1 follows.

2. Statement 2 follows from statement 1:
 \begin{eqnarray*}
 \s (H)&=&\s (\diag (x_1,\ldots , x_n)\der ) =\s (\diag (x_1,\ldots , x_n)) \s (\der )=\diag (\s (x_1),\ldots , \s (x_n))C_\s \der \\
  &=& \diag (\s (x_1),\ldots , \s (x_n)) (\diag (\s (x_1)^{-1},  \ldots , \s(x_n)^{-1})A^{-1} \diag (x_1,\ldots , x_n)= A^{-1}H.\;\; \Box
\end{eqnarray*}


Let  group $G$ acts on a set $S$ and $T\subseteq S$. Then $\Fix_G(T):=\{ g\in G\, | \, gt=t$ for all $t\in T\}$ is the {\em fixator} of the the set $T$. $\Fix_G(T)$ is a subgroup of $G$.

\begin{lemma}\label{k21Mar13}
Let $\s \in \mmWn$. Then $\s (H) = A_{\s^{-1}} H$ for some $A_{\s^{-1}} \in \GL_n(\Z )$ (see (\ref{sHAs}) and Lemma \ref{f21Mar13}) and $\s_{A_{\s^{-1}}} \s \in \Fix_{\mmWn}(H_1, \ldots , H_n)$ where $\s_{A_{\s^{-1}}}\in \GL_n(\Z ) \subseteq \mLn$, see (\ref{msAG}).
\end{lemma}

{\it Proof}. The statement follows from Lemma \ref{h21Mar13}.(2).  $\Box $


$$\Sh_n:=\{ s_\l \in G_n\, | \, s_\l (x_1)=x_1+\l_1, \ldots , s_\l (x_n) = x_n+\l_n\}$$ is the {\em shift group} of automorphisms of the polynomial algebra $P_n$ where  $\l = (\l_1, \ldots , \l_n)\in K^n$; $\Sh_n\subset \AutKalg (P_n) \subseteq \AutLie (D_n)$.

\begin{proposition}\label{WB11Mar13}
$\Fix_{\mmWn}(\der_1, \ldots , \der_n ) = \{ e\}$.
\end{proposition}

{\it Proof}.   Let $\s\in F:= \Fix_{\mmWn} (\der_1, \ldots , \der_n)$. We have to show that $\s = e$.
Let $N:= \Nil_{W_n}(\der_1, \ldots , \der_n):= \{ w\in W_n\,  | \, \d_i^s(w)=0$ for some $s=s(w)$ and all $i=1, \ldots , n\}$. Clearly, $N= D_n$. The automorphisms $\s$ and $\s^{-1}$  preserve the space $N=D_n$, that is $\s^{\pm 1} (D_n) \subseteq D_n$. Hence $\s (D_n) = D_n$ and $\s |_{D_n} \in \Fix_{G_n} (\der_1, \ldots , \der_n)= \Sh_n$, \cite{Bav-Aut-Der-Pol}. The only element $s_\l$ of $\Sh_n$ that can be extended to an automorphism of $W_n$ is $e$ (since $s_\l (x_i^{-1}\der_i) = (x_i+\l_1)^{-1}\der_i$). Therefore,  $\s =e$.
 In more detail, suppose that $s_\l$ can be extended to an automorphism of the Witt algebra $W_n$ and $\l_i\neq 0$, we seek a contradiction. Then applying $s_\l$ to the relation $[x_i^{-1}\der_i , x_i^2\der_i]=3\der_i$ we obtain the relation $ [s_\l (x_i^{-1}\der_i) , (x_i+\l_i)^2\der_i]=3\der_i$. On the other hand, $[(x_i+\l_i)^{-1}\der_i , (x_i+\l_i)^2\der_i]=3\der_i$  in the Lie algebra $K(x_i)\der_i$. Hence, $s_\l (x_i^{-1}\der_i) - (x_i+\l_i)^{-1}\der_i\in C:= C_{K(x_i)\der_i}((x_i+\l_i)^2\der_i)$. Since $C= K\cdot (x+\l_i)^2\der_i$, we see that $s_\l (x_i^{-1}\der_i) \not\in W_n$, a contradiction. In more detail, let $\alpha =(x_i+\l_i)^2$. Then $\beta \der_i\in C$ where $\beta \in K(x_i)$ iff (where $\alpha':= \frac{d \alpha }{dx_i}$, etc) $0=[\alpha \der_i , \beta \der_i]=(\alpha\beta'-\alpha' \beta ) \der_i = \alpha^2(\frac{\beta }{\alpha })'\der_i$ iff $(\frac{\beta }{\alpha })'=0$ iff $\frac{\beta }{\alpha }\in \ker_{K(x_i)}(\der_i) = K$. Hence, $\beta \in K \alpha$, as required.   $\Box $


\begin{lemma}\label{l21Mar13}
$\Fix_{\mmWn} (H_1, \ldots , H_n)=\mT^n$.
\end{lemma}

{\it Proof}. The inclusion $\mT^n \subseteq F:=\Fix_{\mmWn} (H_1, \ldots , H_n)$ is obvious. Let $\s\in F$. We have to show that $\s \in \mT^n$. In view of Proposition \ref{WB11Mar13}, it suffices to show that $\s (\der_1) = \l_1\der_1, \ldots , \s (\der_n) = \l_n\der_n$ for some $\l = (\l_1, \ldots , \l_n) \in K^{*n}$ since then $t_\l\s \in \Fix_{\mmWn} (\der_1, \ldots , \der_n) = \{ e\}$ (Proposition \ref{WB11Mar13}), and so $ \s = t_\l^{-1} \in \mT^n$. Since $\s \in F$, the automorphism respects the weight components of the Lie algebra $W_n$, that is $\s (x^\alpha \CH_n) = x^\alpha \CH_n$ for all $\alpha \in \Z^n$. In particular, for $i=1, \ldots , n$,
\begin{equation}\label{xHss}
 \der_i' = \s (\der_i)= \s ( x_i^{-1} H_i) = x_i^{-1} \sum_{j=1}^n \l_{ij} H_j= -x_i^{-1} \sum_{j=1}^n \l_{ij} x_j\der_j,
\end{equation}
$ \der'= D^{-1}\L D\der$ where $D=\diag (x_1, \ldots ,x_n)$ and $D^{-1}\L D\in \GL (L_n)$, and so  $\L = (\l_{ij})\in \GL_n(K)$. In view of (\ref{xHss}), e have to show that $\L$ is a diagonal matrix. The elements $\der_1, \ldots ,\der_n$ commute, so do $\der_1', \ldots , \der_n'$: for all $i,j=1, \ldots , n$,
$$0=[\der_i', \der_j']= [ x_i^{-1}\sum_{k=1}^n \l_{ik}H_k, x_j^{-1}\sum_{l=1}^n \l_{jl}H_l].$$
Therefore, foro all $i,j, l=1, \ldots , n$,  $\l_{ij}\l_{jl}=\l_{ji}\l_{il}$. For each $i=1, \ldots , n$, let $c_i:= \sum_{j=1}^n \l_{ji}$. The above equalities yield the equalities
$$ \sum_{j=1}^n \l_{ij}\l_{jl} = c_i\l_{il}\;\; {\rm for}\;\; i,l=1, \ldots , n.$$Equivalently, $\L^2= \diag (c_1, \ldots , c_n)\L$. Therefore, $\L =\diag (c_1, \ldots , c_n)$ since $\L \in \GL_n(K)$, as required.  $\Box $

$\noindent $

{\bf Proof of Theorem \ref{20Mar13}}. Let $\s \in \mmWn$. We have to show that $\s \in \mLn$. By Lemma \ref{k21Mar13} and Lemma \ref{l21Mar13}, $\tau \s \in \Fix_{\mmWn} ( H_1, \ldots , H_n)= \mT^n$ for some $\tau \in \mLn$, hence $\s \in \mLn$. $\Box$


\section{The group of automorphisms of the Virasoro algebra  }\label{AUT-Vir}

The aim of this section is to find the group of automorphisms of the Virasoro algebra (Theorem \ref{22Mar13}). The key idea is to use Theorem \ref{22Mar13}.

$\noindent $

{\bf Proof of Theorem \ref{22Mar13}}. Let a $K$-linear map $s: W\ra \CG$ be a
 section to  the surjection $\pi : \CG \ra W$, $a\mapsto a+Z$, i.e. $\pi s = \id_W$. The map $\s$ is an injection  and we identify the vector space $W$ with its image in $\CG$ via $s$. Then, $\CG = W\oplus Z$, a direct sum of vector spaces.

(i) $\widehat{\s}$ {\em is unique}:  Suppose we have another extension, say $\widehat{\s}_1$. Then $\tau := \widehat{\s}_1^{-1}\widehat{\s}\in G:= \AutLie (\CG )$ and
$$\phi (w ) := \tau (w)-w\in Z \;\; {\rm  for\; all}\;\; w\in W,$$
 where $\phi\in \Hom_K(W, Z)$. By condition 2, the inclusion $Z\subseteq [ \CG , \CG ]=[W+Z, W+Z]=[W,W] $ implies that $\tau (z) = z$ for all $z\in Z$. For all $w_1, w_2\in W$,
\begin{equation}\label{w1w2}
[w_1,w_2] = [w_1,w_2]_W+z(w_1, w_2)
\end{equation}
where $[\cdot, \cdot ]$ and  $[\cdot, \cdot ]_W $ are  the Lie brackets in $\CG$ and $W$ respectively and $z(w_1,w_2)\in Z$. Moreover, $ [w_1,w_2]_W$ means $s ( [w_1,w_2]_W)$. Applying the automorphism $\tau$ to the above equality we have
\begin{eqnarray*}
 [w_1,w_2]&=& [\tau (w_1),\tau (w_2)]= \tau ([w_1,w_2]) = \tau ([w_1,w_2]_W+z(w_1,w_2))\\
 &=& [w_1,w_2]_W+\phi ([w_1,w_2]_W) +z(w_1,w_2)=[w_1,w_2]+\phi ([w_1,w_2]_W).
\end{eqnarray*}
Hence, $\phi ([w_1,w_2]_W)=0$ for all $w_1, w_2\in W$. By condition 3, $\phi =0$, that is $\tau (w) = w$ for all $w\in W$. Together with the condition $\tau (z) = z$ for all $z\in Z$, this gives $\tau = e$. So, $\widehat{\s}=\widehat{\s}_1$.

(ii) {\em The map $\s \mapsto \widehat{\s}$ is a monomorphism}: Let $\widehat{\s}$ and $\widehat{\tau}$ be the extensions of $\s $ and $\tau$ respectively. By the uniqueness, $\widehat{\s}\widehat{\tau}$ is the extension of $\s \tau$, that is $\widehat{\s\tau}=\widehat{\s}\widehat{\tau}$, and so the map $\s \mapsto \widehat{\s}$ is a homomorphism. Again, by the uniqueness, $\s \mapsto \widehat{\s}$ is a monomorphism.

(iii) {\em The map $\s \mapsto \widehat{\s}$ is an isomorphism}: By condition 1, the map $\s \mapsto \widehat{\s}$ is a surjection, hence an isomorphism, by (ii). $\Box$

$\noindent $

{\bf Proof of Theorem \ref{A20Mar13}}.  The conditions of Theorem \ref{22Mar13} are satisfied for the Virasoro algebra: $Z= Z(\Vir ) = Kc$, $\Vir / Z\simeq W_1$, $[W_1, W_1]=W_1$ (since $W_1$ is a simple Lie algebra),  $Z\subseteq [ \Vir, \Vir ] $ and each automorphism $\s \in \mmW1=\AutLie (W_1) = \AutKalg (L_1)= \GL_1(\Z ) \ltimes \mT^1$ is extended to an automorphism $\widehat{\s}\in \AutLie (\Vir )$ by the rule $\widehat{\s} (c)=c$. The last condition is obvious for $\s \in \mT^1$ but for $e\neq \s \in \GL_1(\Z ) = \{ \pm 1\}$, i.e. $\s : L_1\ra L_1$, $x\mapsto x^{-1}$, i.e. $\s : W_1\ra W_1$, $ x^iH\mapsto - x^{-1}H$ for all $i\in \Z$, it follows from the relations (\ref{defVir}). $\Box $


\begin{corollary}\label{pA20Mar13}

\begin{enumerate}
\item Each automorphism $\s $ of the Witt algebra $W_1$ is uniquely extended to an automorphism $\widehat{\s}$ of the Virasoro algebra $\Vir$. Moreover, $\widehat{\s}(c)= c$.
\item All the automorphisms of the Virasoro algebra $\Vir$ act trivially on its centre.
\end{enumerate}
\end{corollary}

When we drop condition 3 of Theorem \ref{22Mar13}, we obtain a more general result.
\begin{corollary}\label{AA22Mar13}
Let $\CG$ be a Lie algebra, $Z$ be a subspace of the centre of $\CG$  and $W=\CG / Z$.  Fix a $K$-linear map $s: W\ra \CG$ which is a section to the surjection $\pi : \CG \ra W$, $a\mapsto a+Z$, and identify $W$ with $\im (s)$, and so $\CG = W\oplus Z$ (a  direct sum of vector spaces). Let $\CK := \{ \tau = \tau_\phi \in \End_K(\CG ) \, | \, \tau (w) = w+\phi (w)$ and $\tau (z) = z$ for all $w\in W$ and $z\in Z$, $\phi \in \Hom_K( W, Z)$ is such that $\phi ([W, W]) =0\}$. Suppose that
\begin{enumerate}
\item every automorphism $\s$ of the Lie algebra $W$ can be extended to an automorphism $\widehat{\s}$ of the Lie algebra $\CG$, and
\item $Z\subseteq [G, G]$.
\end{enumerate}
Then the short exact sequence of groups
$$ 1\ra \CK \ra \AutLie (\CG ) \stackrel{\psi}{\ra} \AutLie (W)\ra 1$$ is exact where $\psi ( \s ) : a+Z\mapsto \s (a) +Z$ for all   $a\in \CG$.
\end{corollary}

{\it Proof}. By condition 1, $\psi$ is a group epimorphism. It remains to show that $\ker (\psi ) = \CK$. Let $\tau \in \ker (\psi )$. Each element $g\in \CG = W\oplus Z$ is a unique sum $g=w+z$ where $w\in W$ and $z\in Z$. Then $\tau (w) = w+\phi (w)$ for some $\phi\in\Hom_K(W, Z)$. We keep the notation of the proof of Theorem \ref{22Mar13}. By condition 2, $Z\subseteq [ \CG , \CG ] = [W,W]$, hence $ \tau (z) = z$ for all elements $z\in Z$. Applying the automorphism $\tau $ to the equality (\ref{w1w2}) yields $\phi ([w_1, w_2]) =0$ (see the proof of Theorem \ref{22Mar13}). It follows that $\ker (\psi ) = \CK$.  $\Box $



$${\bf Acknowledgements}$$

 The work is partly supported by  the Royal Society  and EPSRC.

\small{

Department of Pure Mathematics

University of Sheffield

Hicks Building

Sheffield S3 7RH

UK

email: v.bavula@sheffield.ac.uk}

\end{document}